\documentclass[preprint,12pt]{elsarticle}

\usepackage{amssymb}
\usepackage{amsthm}
\usepackage{graphics}
\usepackage{epsfig}
\usepackage{graphicx}
\usepackage{color}
\usepackage{hyperref}
\usepackage{a4wide}
\usepackage{amsmath}
\usepackage{amsfonts}
\usepackage{enumerate}
\newcommand{\pf}{\noindent {\bf Proof: }}

\newtheorem{lemma}{\bf Lemma}[section]
 \newtheorem{example}{\bf Example}[section]
\newtheorem{theorem}{\bf Theorem}[section]

\newtheorem{corollary}[lemma]{\bf Corollary}

\newtheorem{remark}{\bf Remark}[section]

\journal{Communications in Algebra}

\begin{document}

\begin{frontmatter}



\title{Subspace Sum Graph of a Vector Space}



\author{Angsuman Das\corref{cor1}}
\ead{angsumandas@sxccal.edu}

\address{Department of Mathematics,\\ St. Xavier's College, Kolkata, India.\\angsumandas@sxccal.edu}
\cortext[cor1]{Corresponding author}


\begin{abstract}
In this paper we introduce a graph structure, called subspace sum graph $\mathcal{G}(\mathbb{V})$ on a finite dimensional vector space $\mathbb{V}$ where the vertex set is the collection of non-trivial proper subspaces of a vector space and two vertices $W_1,W_2$ are adjacent if $W_1 + W_2=\mathbb{V}$. The diameter, girth, connectivity, maximal independent sets, different variants of domination number, clique number and chromatic number of $\mathcal{G}(\mathbb{V})$ are studied. It is shown that two subspace sum graphs are isomorphic if and only if the base vector spaces are isomorphic. Finally some properties of subspace sum graph are studied when the base field is finite.
\end{abstract}


\end{frontmatter}


\section{Introduction}
Apart from its combinatorial motivation, graph theory also helps to characterize various algebraic structures by means of studying certain graphs associated to them. Till date, a lot of research, e.g., \cite{anderson-livingston,mks-ideal,survey2} has been done in connecting graph structures to various algebraic objects. Recently, some works associating graphs with subspaces of vector spaces can be found in \cite{angsu-comm-alg,angsu-lin-mult-alg,angsu-comm-alg-2,angsu-jaa,int-vecsp-2,int-vecsp-3,int-vecsp-1}. 

In this paper we define a graph structure on a finite dimensional vector space $\mathbb{V}$ over a field $\mathbb{F}$, called Subspace Sum Graph of $\mathbb{V}$ and derive some properties of the graph using the algebraic properties of vector subspaces.

\section{Definitions and Preliminaries}
In this section, for convenience of the reader and also for later use, we recall some definitions, notations and results concerning elementary graph theory. For undefined terms and concepts the reader is referred to \cite{west-graph-book}.

By a graph $G=(V,E)$, we mean a non-empty set $V$ and a symmetric binary relation (possibly empty) $E$ on $V$. The set $V$ is called the set of vertices and $E$ is called the set of edges of $G$. Two element $u$ and $v$ in $V$ are said to be adjacent if $(u,v) \in E$. $H=(W,F)$ is called a {\it subgraph} of $G$ if $H$ itself is a graph and $\phi \neq W \subseteq V$ and $F \subseteq E$. If $V$ is finite, the graph $G$ is said to be finite, otherwise it is infinite. Two graphs $G=(V,E)$ and $G'=(V',E')$ are said to be {\it isomorphic} if there exists a bijection $\phi: V \rightarrow V'$ such that $(u,v) \in E \mbox{ iff } (\phi(u),\phi(v)) \in E'$. A {\it path} of length $k$ in a graph is an alternating sequence of vertices and edges, $v_0,e_0,v_1,e_1,v_2,\ldots, v_{k-1},e_{k-1},v_k$, where $v_i$'s are distinct (except possibly the first and last vertices) and $e_i$ is the edge joining $v_i$ and $v_{i+1}$. We call this a path joining $v_0$ and $v_{k}$. A {\it cycle} is a path with $v_0=v_k$. A graph is said to be {\it Eulerian} if it contains a cycle containing all the edges in $G$ exactly once. A cycle of length 3 is called a {\it triangle}. A graph is {\it connected} if for any pair of vertices $u,v \in V,$ there exists a path joining $u$ and $v$. A graph is said to be {\it triangulated} if for any vertex $u$ in $V$, there exist $v,w$ in $V$, such that $(u,v,w)$ is a triangle. The {\it distance} between two vertices $u,v \in V,~ d(u,v)$ is defined as the length of the shortest path joining $u$ and $v$, if it exists. Otherwise, $d(u,v)$ is defined as $\infty$. The {\it diameter} of a graph is defined as $diam(G)=\max_{u,v \in V}~ d(u,v)$, the largest distance between pairs of vertices of the graph, if it exists. Otherwise, $diam(G)$ is defined as $\infty$. The {\it girth} of a graph is the length of its shortest cycle, if it exists. Otherwise, it is defined as $\infty$. If all the vertices of $G$ are pairwise adjacent, then $G$ is said to be {\it complete}. A complete subgraph of a graph $G$ is called a {\it clique}. A {\it maximal clique} is a clique which is maximal with respect to inclusion. The {\it clique number} of $G$, written as $\omega(G)$, is the maximum size of a clique in $G$. A subset $I$ of $V$ is said to be {\it independent} if any two vertices in that subset are pairwise non-adjacent. A {\it maximal independent set} is an independent set which is maximal with respect to inclusion. The {\it chromatic number} of $G$, denoted as $\chi(G)$, is the minimum number of colours needed to label the vertices so that the adjacent vertices receive different colours. It is known that for any graph $G$, $\chi(G)\geq \omega(G)$.  A graph $G$ is said to be {\it weakly perfect} if $\chi(G)= \omega(G)$. A graph $G$ is said to be {\it perfect} if $\chi(H)= \omega(H)$, for all induced subgraph $H$ of $G$.  A subset $D$ of $V$ is said to be {\it dominating set} if any vector in $V \setminus D$ is adjacent to at least one vertex in $D$. A subset $D$ of $V$ is said to be {\it total dominating set} if any vector in $V$ is adjacent to at least one vertex in $D$. A dominating set $D$  of $G$ is said to be a connected dominating set of $G$ if the subgraph generated by $D$, $\langle D \rangle$ is connected. A dominating set $D$  of $G$ is said to be a dominating clique of $G$ if the subgraph generated by $D$, $\langle D \rangle$ is complete. The {\it dominating number} $\gamma(G)$, the {\it total dominating number} $\gamma_t(G)$, the {\it connected dominating number} $\gamma_c(G)$ and the {\it clique domination number} are the minimum size of a dominating set, a total dominating set, a connected dominating set and a dominating clique in $G$ respectively.

\section{Subspace Sum Graph of a Vector Space}
Let $\mathbb{V}$ be a finite dimensional vector space over a field $\mathbb{F}$ of dimension greater than $1$ and $\theta$ denote the null vector. We define a graph $\mathcal{G}(\mathbb{V})=(V,E)$ as follows: $V=$ the collection of non-trivial proper subspaces of $\mathbb{V}$ and for $W_1,W_2 \in V$, $W_1 \sim W_2$ or $(W_1,W_2) \in E$ if $W_1 + W_2=\mathbb{V}$. Since, $dim(\mathbb{V})>1$, $V \neq \emptyset$.

\begin{figure}[ht]
\centering
\begin{center}
\includegraphics[scale=.4]{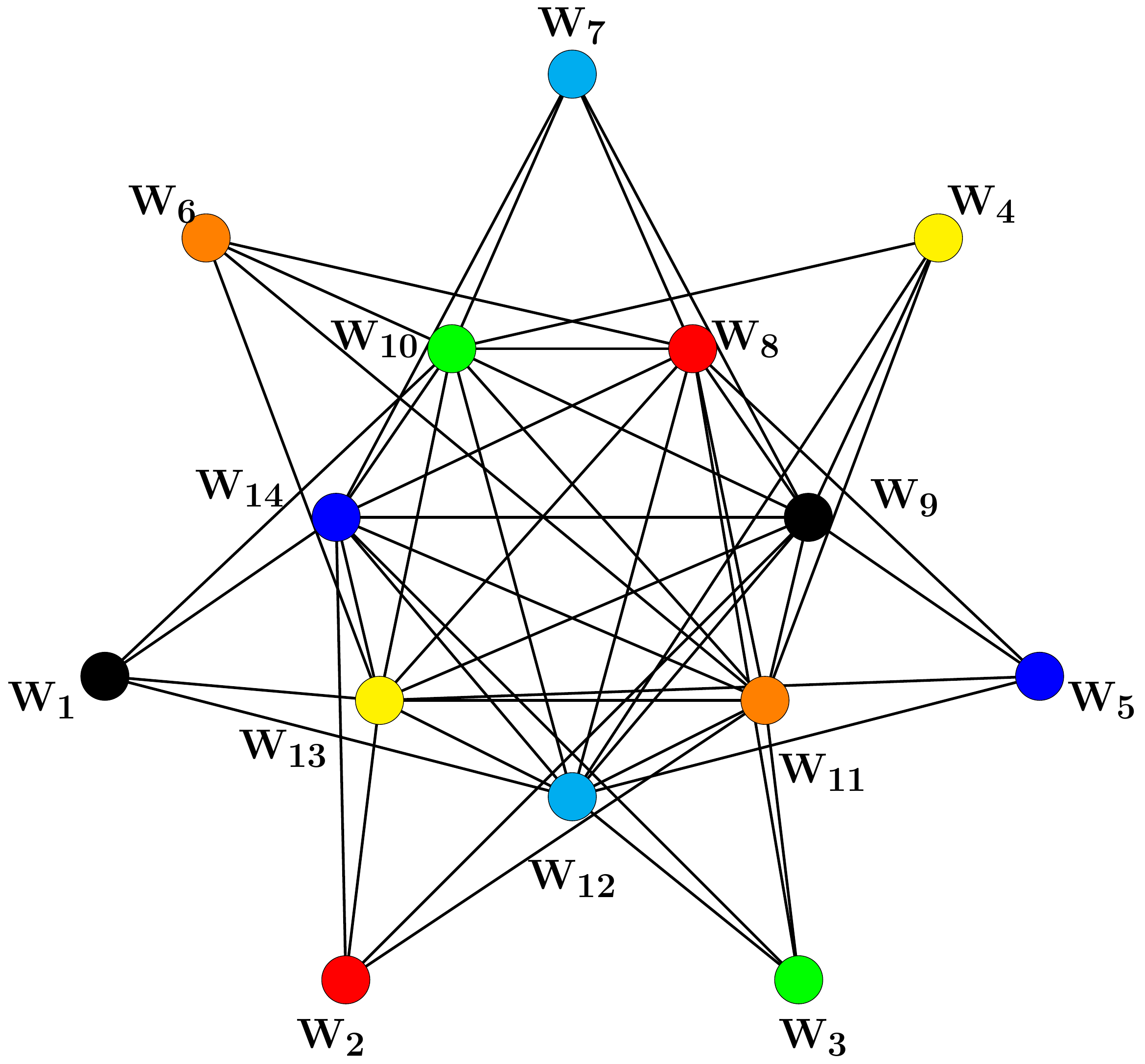}
\caption{$\mathcal{G}(\mathbb{V})$ with $dim(\mathbb{V})=3; \mathbb{F}=\mathbb{Z}_2$}
\label{diagram}
\end{center}
\end{figure}

{\example \label{example} Consider a $3$ dimensional vector space $\mathbb{V}$ over $\mathbb{Z}_2$ with a basis $\{\alpha_1,\alpha_2,\alpha_3\}$. The following are the possible non-trivial proper subspaces of $\mathbb{V}$: $W_1=\langle \alpha_1 \rangle,W_2=\langle \alpha_2 \rangle$, $W_3=\langle \alpha_3 \rangle,W_4=\langle \alpha_1 +\alpha_2 \rangle,W_5=\langle \alpha_2+\alpha_3 \rangle,W_6=\langle \alpha_1+\alpha_3 \rangle,W_7=\langle \alpha_1+\alpha_2+\alpha_3 \rangle,W_8=\langle \alpha_1,\alpha_2 \rangle,W_9=\langle \alpha_1,\alpha_3 \rangle,W_{10}=\langle \alpha_2,\alpha_3 \rangle,W_{11}=\langle \alpha_1,\alpha_2+\alpha_3 \rangle,W_{12}=\langle \alpha_2, \alpha_1+\alpha_3 \rangle,W_{13}=\langle \alpha_3,\alpha_1+\alpha_2 \rangle,W_{14}=\langle \alpha_1+\alpha_2,\alpha_2+\alpha_3 \rangle$. Then the graph $\mathcal{G}(\mathbb{V})$ is given in Figure \ref{diagram}. For discussion on chromatic number and clique number of this graph, see Remark \ref{example-remark} of Section \ref{finite-field}. }

Throughout this paper, even if it is not mentioned explicitly, the underlying field is $\mathbb{F}$ and $\mathbb{V}$ is finite dimensional. Now we study some basic properties like completeness, connectedness, diameter and girth of $\mathcal{G}(\mathbb{V})$.

{\theorem \label{complete} $\mathcal{G}(\mathbb{V})$ is complete if and only if $dim(\mathbb{V})=2$.}\\
\pf If $dim(\mathbb{V})=2$, then the vertices of $\mathcal{G}(\mathbb{V})$ are the one dimensional subspaces of $\mathbb{V}$. Now, the sum of two distinct one dimensional subspaces in a two dimensional vector space is two dimensional and hence equal to $\mathbb{V}$ and hence $\mathcal{G}(\mathbb{V})$ is complete. Conversely, if $\mathcal{G}(\mathbb{V})$ is complete and $dim(\mathbb{V})\geq 2$, then there exists two distinct one dimensional subspaces of $\mathbb{V}$ whose sum is not $\mathbb{V}$. This contradicts the completeness of $\mathcal{G}(\mathbb{V})$ and hence $dim(\mathbb{V})=2$. \qed

{\theorem \label{diameter} If $dim(\mathbb{V})\geq 3$, then $\mathcal{G}(\mathbb{V})$ is connected and  $diam(\mathcal{G}(\mathbb{V}))= 2$.}\\
\pf Let $W_1, W_2$ be two distinct non-trivial proper subspaces of $\mathbb{V}$. If $W_1 + W_2 = \mathbb{V}$, then $d(W_1, W_2)=1$ and we are done. However, as $dim(\mathbb{V})\geq 3$, by Theorem \ref{complete}, there exists $W_1, W_2$ which are not adjacent in $\mathcal{G}(\mathbb{V})$. If $W_1 \not\sim W_2$, then $W_1 + W_2 \subset \mathbb{V}$. Two cases may arise.\\
 {\bf Case 1:} $W_1 \cap W_2 \neq \{\theta \}$. Then there exists a non-trivial proper subspace $W_3$ such that $(W_1 \cap W_2)+W_3=\mathbb{V}$. Since, $W_1 \cap W_2 \subset W_1,W_2$, therefore $W_1 + W_3=\mathbb{V}$ and $W_2+ W_3 = \mathbb{V}$. Thus, we have $W_1 \sim W_3 \sim W_2$ and hence $d(W_1,W_2)=2$.\\
{\bf Case 2:} $W_1 \cap W_2 = \{\theta \}$. Let $dim(W_1)=k_1; dim(W_2)=k_2$ with $k_1 \leq k_2$. Moreover, let $W_1=\langle \alpha_1,\alpha_2,\ldots,\alpha_{k_1} \rangle$; $W_2=\langle \beta_1,\beta_2,\ldots,\beta_{k_2} \rangle$. Since $W_1 \cap W_2 = \{\theta \}$, therefore $\{ \alpha_1,\alpha_2,\ldots,\alpha_{k_1},\beta_1,\beta_2,\ldots,\beta_{k_2} \}$ is linearly independent. Also, as $W_1 + W_2 \subset \mathbb{V}$, $\{ \alpha_1,\alpha_2,\ldots,\alpha_{k_1},\beta_1,\beta_2,\ldots,\beta_{k_2} \}$ can be extended to a basis $$\{ \alpha_1,\alpha_2,\ldots,\alpha_{k_1},\beta_1,\beta_2,\ldots,\beta_{k_2},\gamma_1,\gamma_2,\ldots,\gamma_{k_3} \}\mbox{ of }\mathbb{V}, \mbox{ where }dim(\mathbb{V})=n=k_1 + k_2 + k_3.$$ Set $W_3=\langle \alpha_1 + \beta_1, \alpha_2 + \beta_2,\ldots, \alpha_{k_1}+\beta_{k_1},\beta_{k_1 + 1},\beta_{k_1 + 2},\ldots,\beta_{k_2},\gamma_1,\gamma_2,\ldots,\gamma_{k_3} \rangle$. Clearly, $dim(W_3)=k_2+k_3<n$. We claim that $W_1+W_3=\mathbb{V}$ and $W_2+W_3=\mathbb{V}$.\\
{\bf Proof of Claim:} $\beta_1=(\alpha_1+\beta_1)-\alpha_1 \in W_1+W_3$. Similarly, $\beta_2,\beta_3,\ldots,\beta_{k_1} \in W_1+W_3$. Also, $\beta_{k_1 +1},\ldots, \beta_{k_2},\gamma_1,\gamma_2,\ldots,\gamma_{k_3} \in W_3$ and $\alpha_1,\alpha_2,\ldots,\alpha_{k_1} \in W_1$. Thus,\\ 
$\langle \alpha_1,\alpha_2,\ldots,\alpha_{k_1},\beta_1,\beta_2,\ldots,\beta_{k_2},\gamma_1,\gamma_2,\ldots,\gamma_{k_3} \rangle \subset W_1+W_3$ and hence $W_1+W_3=\mathbb{V}$. Similarly $W_2+W_3=\mathbb{V}$.

Thus we have $W_1 \sim W_3 \sim W_2$, and combining both the cases, we get  $diam(\mathcal{G}(\mathbb{V}))= 2$.\qed

{\theorem $\mathcal{G}(\mathbb{V})$ is triangulated and hence $girth(\mathcal{G}(\mathbb{V}))=3$.}\\
\pf Let $W_1=\langle \alpha_1,\alpha_2,\ldots,\alpha_k \rangle$ be a $k$-dimensional subspace of $\mathbb{V}$ and $\{\alpha_1,\alpha_2,\ldots,\alpha_k\}$ be extended to a basis of $\mathbb{V}$ as $\{\alpha_1,\alpha_2,\ldots,\alpha_k,\beta_1,\beta_2,\ldots,\beta_l \}$  where $n=k+l$. Without loss of generality, we take $k \leq l$ and set $W_2=\langle \beta_1, \beta_2, \ldots, \beta_l \rangle$ and $W_3= \langle \alpha_1 + \beta_1, \alpha_2+ \beta_2, \ldots, \alpha_k + \beta_k, \beta_{k+1},\ldots, \beta_l \rangle$. By similar arguments to that of the proof of claim in Theorem \ref{diameter}, it is clear that $W_1+W_3=W_2+W_3=W_1+W_2=\mathbb{V}$ and hence $W_1,W_2,W_3$ forms a triangle in $\mathcal{G}(\mathbb{V})$. \qed

\section{Maximal Cliques and Maximal Independent Sets in $\mathcal{G}(\mathbb{V})$}
In this section, we study the structure of maximal cliques and maximal independent sets in $\mathcal{G}(\mathbb{V})$.

{\theorem \label{maximal-clique} The collection $\mathcal{W}_{n-1}$ of all $(n-1)$ dimensional subspaces of $\mathbb{V}$ is a maximal clique in $\mathcal{G}(\mathbb{V})$.}\\
\pf As any two distinct $n-1$ dimensional subspaces are adjacent in $\mathcal{G}(\mathbb{V})$, the collection $\mathcal{W}_{n-1}$ of all $(n-1)$ dimensional subspaces of $\mathbb{V}$ is a clique in $\mathcal{G}(\mathbb{V})$. For maximality, let if possible, $\mathcal{W}_{n-1} \cup \{W\}$ be a clique where $W$ is non-trivial proper subspace of $\mathbb{V}$ with $dim(W)=k < n-1$. Since $k < n-1$, there exist $W' \in \mathcal{W}_{n-1}$ such that $W \subset W'$. Thus $W+W'=W'\neq \mathbb{V}$, thereby contradicting that $\mathcal{W}_{n-1} \cup \{W\}$ is a clique. Hence, $\mathcal{W}_{n-1}$ is a maximal clique in $\mathcal{G}(\mathbb{V})$. \qed 

{\lemma \label{graph-homomorphism-lemma} There exists a graph homomorphism from $\mathcal{G}(\mathbb{V})$ to the subgraph induced by $\mathcal{W}_{n-1}$, defined in Theorem \ref{maximal-clique}. }\\
\pf Since $\mathcal{W}_{n-1}$ is a clique, subgraph induced by it $\langle \mathcal{W}_{n-1} \rangle$ is a complete graph. As every non-trivial proper subspace $W$ of $\mathbb{V}$ is contained in at least one (possibly more than one) subspace in $\mathcal{W}_{n-1}$,  there exists a map $\varphi: \mathcal{G}(\mathbb{V}) \rightarrow \langle \mathcal{W}_{n-1}\rangle$ given by $W \mapsto W'$, where $W'$ is an $n-1$ dimensional subspace of $\mathbb{V}$ containing $W$. The existence of such a map is guaranteed by the axiom of choice. 

It is to be noted that if $W_1 \sim W_2$ in $\mathcal{G}(\mathbb{V})$, i.e., $W_1+W_2=\mathbb{V}$, then $W_1$ and $W_2$ are not contained in same subspace in $\mathcal{W}_{n-1}$. Thus $W'_1\neq W'_2$. Thus, as $\langle \mathcal{W}_{n-1} \rangle$ is a complete graph, $W'_1\sim W'_2$. Hence $\varphi$ preserve adjacency and is a graph homomorphism. \qed

{\remark Lemma \ref{graph-homomorphism-lemma} will be useful in finding clique number and chromatic number of $\mathcal{G}(\mathbb{V})$ when the field $\mathbb{F}$ is finite. (See Theorem \ref{clique-number-chromatic-number-theorem})}

{\theorem \label{maximal-ind-set-odd-case} If $n$ is odd, i.e., $n=2m+1$, then the collection $\mathcal{W}[{m}]$ of all non-trivial subspaces of $\mathbb{V}$ with dimension less than or equal to $m$ is a maximal independent set in $\mathcal{G}(\mathbb{V})$.}\\
\pf For any $W_1,W_2 \in \mathcal{W}[{m}]$, $W_1+W_2 \neq \mathbb{V}$ as $dim(W_1+W_2)\leq 2m <n$. Thus, $\mathcal{W}[{m}]$ is an independent set. For any $W$ with $dim(W)>m$, there exists a subspace $W' \in \mathcal{W}[{m}]$ such that $W+W'=\mathbb{V}$. Thus, $\mathcal{W}[{m}]$ is a maximal independent set. \qed

{\theorem \label{maximal-ind-set-even-case} Let $n$ is even, i.e., $n=2m$ and $\mathcal{W}[{m-1}]$ be the collection of all non-trivial subspaces of $\mathbb{V}$ with dimension less than or equal to $m-1$. Let $\alpha (\neq \theta) \in \mathbb{V}$ and $\mathcal{W}^\alpha[{m}]$ be the collection of all non-trivial subspaces of $\mathbb{V}$ containing $\alpha$ and having dimension $m$. Then $\mathcal{W}[{m-1}]\cup \mathcal{W}^\alpha[{m}]$ is a maximal independent set in $\mathcal{G}(\mathbb{V})$.}\\
\pf Let $W_1,W_2 \in \mathcal{W}[{m-1}]\cup \mathcal{W}^\alpha[{m}]$. If at least one of $W_1,W_2 \in \mathcal{W}[{m-1}]$, $W_1+W_2 \neq \mathbb{V}$ as $dim(W_1+W_2)\leq 2m-1 <n$ and thus $W_1 \not\sim W_2$. If both $W_1,W_2 \in \mathcal{W}^\alpha[{m}]$, then as $\alpha \in W_1 \cap W_2$, $dim(W_1\cap W_2)\geq 1$ and as a result $dim(W_1+W_2)=dim(W_1)+dim(W_2)-dim(W_1\cap W_2)\leq 2m-1<n$. Thus $W_1 \not\sim W_2$. Hence $\mathcal{W}[{m-1}]\cup \mathcal{W}^\alpha[{m}]$ is an independent set.

For maximality, let $W$ be a proper subspace of $\mathbb{V}$ not in $\mathcal{W}[{m-1}]\cup \mathcal{W}^\alpha[{m}]$. If $dim(W)=m$, then as $\alpha \not\in W$, there exists a subspace $W' \in \mathcal{W}^\alpha[{m}]$ such that $W+W'=\mathbb{V}$. If $dim(W)>m$, then there exists a subspace $W'' \in \mathcal{W}[{m-1}]$ such that $W+W''=\mathbb{V}$. Thus, $\mathcal{W}[{m-1}]\cup \mathcal{W}^\alpha[{m}]$ is a maximal independent set in $\mathcal{G}(\mathbb{V})$. \qed

\section{Dominating Sets in $\mathcal{G}(\mathbb{V})$}
In this section, we study the minimal dominating sets and domination number of $\mathcal{G}(\mathbb{V})$ and use it prove that the subspace sum graph of two vector spaces are isomorphic if and only if the two vector spaces are isomorphic.
{\lemma \label{dom-set-atleast-n} If $\mathcal{W}$ be a dominating set in $\mathcal{G}(\mathbb{V})$, then $|\mathcal{W}|\geq n$, where $dim(\mathbb{V})=n$.}\\
\pf Let $\alpha$ be a non-null vector in $\mathbb{V}$. Consider the subspace $\langle \alpha \rangle$. Since $\mathcal{W}$ is a dominating set, either $\langle \alpha \rangle \in \mathcal{W}$ or there exists $W_1 \in \mathcal{W}$ such that $\langle \alpha \rangle \sim W_1$ i.e., $\langle \alpha \rangle + W_1 = \mathbb{V}$. If $\langle \alpha \rangle \in \mathcal{W}$, we choose some other $\alpha'\neq \theta$ from $\mathbb{V}$ such that $\langle \alpha' \rangle \not\in \mathcal{W}$. However, if $\langle \alpha \rangle \in \mathcal{W}$ for all $\alpha \in \mathbb{V}\setminus \{\theta\}$, then $\mathcal{W}$ contains all the one-dimensional subspaces of $\mathbb{V}$. Now, irrespective of whether the field is finite or infinite, the number of one-dimensional subspaces of $\mathbb{V}$ is greater than $n$. (See Remark \ref{one-dimensional-count}) Hence, the lemma follows trivially. Thus, we assume that $\alpha \not\in \mathcal{W}$. Then as indicated earlier there exists $W_1 \in \mathcal{W}$ such that $\langle \alpha \rangle + W_1 = \mathbb{V}$. Note that $dim(W_1)=n-1$. 

We choose a non-null vector $\alpha_1 \in W_1$. By similar arguments we can assume that $\langle\alpha_1\rangle \not\in \mathcal{W}$. As $\mathcal{W}$ is a dominating set, there exists $W_2\in \mathcal{W}$ such that $\langle \alpha_1 \rangle + W_2 =\mathbb{V}$. It is to be noted that $dim(W_2)=n-1$ and $W_1\neq W_2$, as $\alpha_1\in W_1\setminus W_2$.

Now consider the subspace $W_1\cap W_2$. As $W_1,W_2$ are of dimension $n-1$, $W_1+W_2=\mathbb{V}$. Thus by the result $dim(W_1+W_2)=dim(W_1)+dim(W_2)-dim(W_1 \cap W_2)$, we have $dim(W_1 \cap W_2)=n-2$. We choose $\alpha_2(\neq \theta) \in W_1 \cap W_2$ and without loss of generality assume that $\langle \alpha_2 \rangle \not\in \mathcal{W}$. Therefore there exists $W_3 \in \mathcal{W}$ such that $\langle \alpha_2 \rangle + W_3=\mathbb{V}$. Again we note that $dim(W_3)=n-1$ and $W_3\neq W_1, W_3 \neq W_2$, as $\alpha_2\in (W_1 \cap W_2)\setminus W_3$.

Now consider the subspace $W_1\cap W_2 \cap W_3$. Since $dim(W_3)=n-1$ and $(W_1 \cap W_2)\not\subset W_3$, we have $(W_1\cap W_2)+W_3=\mathbb{V}$. Thus by using $dim[(W_1\cap W_2)+W_3]=dim(W_1\cap W_2)+dim(W_3)-dim(W_1 \cap W_2 \cap W_3)$, we have $dim(W_1 \cap W_2\cap W_3)=n-3$. As in earlier cases, we choose $\alpha_3(\neq \theta) \in W_1 \cap W_2\cap W_3$ and there exists $W_4 \in \mathcal{W}$ such that $\langle \alpha_3 \rangle + W_4=\mathbb{V}$ with $dim(W_4)=n-1$ and $W_4\neq W_1, W_2, W_3$.

This process continues till we get $W_n \in \mathcal{W}$ such that $\langle \alpha_{n-1} \rangle + W_n=\mathbb{V}$ with $dim(W_n)=n-1$ and $W_n\neq W_1, W_2,\ldots,W_{n-1}$ and $dim(W_1 \cap W_2 \cap \cdots W_n)=0$.

Thus, there are at least $n$ many distinct subspaces $W_1, W_2,\ldots,W_n$ in $\mathcal{W}$ and the lemma holds. \qed

{\remark \label{one-dimensional-count} If the base field $\mathbb{F}$ is infinite and since $n>1$, then there are infinitely many one-dimensional subspaces. If $\mathbb{F}$ is finite with $q$ elements and since $n>1$, then from the proof of Theorem \ref{galois-number-count}, it follows that number of one-dimensional subspaces is $q^{n-1}+q^{n-2}+\cdots+1$ is greater than $n$.}

{\theorem \label{dom-number} Let $\{\alpha_1,\alpha_2,\ldots,\alpha_n\}$ be a basis of $\mathbb{V}$ and $W_i=\langle \alpha_1,\alpha_2,\ldots,\alpha_{i-1},\alpha_{i+1},\ldots \alpha_n\rangle $ for $i=1,2,\cdots,n$. Then the collection $\mathcal{W}=\{W_i: 1\leq i \leq n\}$ is a minimum dominating set in $\mathcal{G}(\mathbb{V})$ and hence $\gamma(\mathcal{G}(\mathbb{V}))=n$.}\\
\pf Let $W$ be an arbitrary vertex of $\mathcal{G}(\mathbb{V})$ and $\theta\neq \alpha \in W$. Suppose $\alpha=c_1\alpha_1+c_2\alpha_2+\cdots +c_n \alpha_n$ with at least one $c_i\neq 0$. If $c_j \neq 0$, then $\langle \alpha \rangle + W_j = \mathbb{V}$ and hence $W+W_j=\mathbb{V}$, i.e., $W \sim W_j$. Thus $\mathcal{W}$ is a dominating set. Now consider $\mathcal{W} \setminus \{W_i\}$. As $\langle \alpha_i \rangle + W_j=W_j$ for $i \neq j$, we have $\langle \alpha_i \rangle \not\sim W_j, \forall W_j \in \mathcal{W}\setminus \{W_i\}$. Thus $\mathcal{W} \setminus \{W_i\}$ is not a dominating set, thereby showing that $\mathcal{W}$ is a minimal dominating set. Now the theorem follows from Lemma \ref{dom-set-atleast-n}. \qed
 
{\corollary Let $\mathbb{V}_1$ and $\mathbb{V}_2$ be two finite dimensional vector spaces over the same field $\mathbb{F}$. Then $\mathbb{V}_1$ and $\mathbb{V}_2$ are isomorphic as vector spaces if and only if $\mathcal{G}(\mathbb{V}_1)$ and $\mathcal{G}(\mathbb{V}_2)$ are isomorphic as graphs.}\\
\pf It is quite obvious that if $\mathbb{V}_1$ and $\mathbb{V}_2$ are isomorphic as vector spaces, then $\mathcal{G}(\mathbb{V}_1)$ and $\mathcal{G}(\mathbb{V}_2)$ are isomorphic as graphs. For the other part, let $\mathcal{G}(\mathbb{V}_1)$ and $\mathcal{G}(\mathbb{V}_2)$ be isomorphic as graphs. Let $dim(\mathbb{V}_1)=n_1$ and $dim(\mathbb{V}_2)=n_2$. Then, by Theorem \ref{dom-number}, the domination numbers of $\mathcal{G}(\mathbb{V}_1)$ and $\mathcal{G}(\mathbb{V}_2)$ are $n_1$ and $n_2$ respectively. However, as the two graphs are isomorphic, $n_1=n_2=n$ (say). Thus $\mathbb{V}_1$ and $\mathbb{V}_2$ are of same finite dimension over the field $\mathbb{F}$ and hence both are isomorphic to $\mathbb{F}^n$ as vector spaces.\qed 

{\corollary If $\mathbb{V}$ is an $n$ dimensional vector space, then the total domination number $\gamma_t$, connected domination number $\gamma_c$ and clique domination number $\gamma_{cl}$ of $\mathcal{G}(\mathbb{V})$ are all equal to $n$, i.e., $$\gamma=\gamma_t=\gamma_c=\gamma_{cl}=n.$$}
\pf It is known that if a graph has a dominating clique and $\gamma \geq 2$, then $\gamma \leq \gamma_t \leq \gamma_c \leq \gamma_{cl}$ (See pg. 167 of \cite{domination-book}). Thus, it suffices to show that the minimum dominating set $\mathcal{W}=\{W_i: 1\leq i \leq n\}$ (constructed as in Theorem \ref{dom-number}) is a total dominating set, a connected dominating set as well as a dominating clique. Since $W_i \sim W_j$ in $\mathcal{G}(\mathbb{V})$ for $i \neq j$, $\mathcal{W}$ is a total dominating set. Moreover, the subgraph $\langle \mathcal{W} \rangle$ spanned by $\mathcal{W}$ is a complete graph and hence connected. Thus, $\mathcal{W}$ is a connected dominating set as well as a dominating clique. \qed

\section{The Case of Finite Fields}\label{finite-field}
In this section, we study some properties of $\mathcal{G}(\mathbb{V})$ if the base field $\mathbb{F}$ is finite, say of order $q=p^r$ where $p$ is a prime. In particular, we find the order, degree, chromatic number, clique number and edge connectivity of $\mathcal{G}(\mathbb{V})$ and show that $\mathcal{G}(\mathbb{V})$ is never Hamiltonian.

It is known that the number of $k$ dimensional subspaces of an $n$-dimensional vector space over a finite field of order $q$ is the $q$-binomial coefficient (See Chapter 7 of \cite{quantum-book})
$$\left[ \begin{array}{c} n \\ k \end{array}\right]_q=\dfrac{(q^n -1)(q^{n}-q)\cdots (q^{n}-q^{k-1})}{(q^k-1)(q^{k}-q)\cdots (q^{k}-q^{k-1})},$$ and hence the total number of non-trivial proper subspaces of $\mathbb{V}$, i.e., the order of $\mathcal{G}(\mathbb{V})$ is given by $$\sum_{k=1}^{n-1}\left[ \begin{array}{c} n \\ k \end{array}\right]_q=G(n,q)-2,$$
where $G(n,q)$ is the Galois number\footnote{For definition, see \cite{galois-number}.}.

{\theorem \label{degree-theorem} Let $W$ be a $k$-dimensional subspace of an $n$-dimensional vector space $\mathbb{V}$ over a finite field $\mathbb{F}$ with $q$ elements. Then degree of $W$ in $\mathcal{G}(\mathbb{V})$, $deg(W)=\sum_{r=0}^{k-1}N_r$, where $$N_r=\dfrac{(q^k-1)(q^k-q)\cdots(q^k-q^{r-1})(q^n-q^k)(q^n-q^{k+1})\cdots(q^n-q^{n-1})}{(q^{n-k+r}-1)(q^{n-k+r}-q)\cdots(q^{n-k+r}-q^{n-k+r-1})}$$}\\
\pf Since $deg(W)=$ the number of subspaces whose sum with $W$ is $\mathbb{V}$ and $dim(W)=k<n$, the subspaces adjacent to $W$ have dimension at least $n-k$, i.e., a subspace $W'$ is adjacent to $W$ if $dim(W')=n-k+r$ and $dim(W\cap W')=r$ where $0\leq r \leq k-1$. To find such subspaces $W'$, we choose $r$ linearly independent vectors from $W$ and $n-k$ linearly independent vectors from $\mathbb{V}\setminus W$, and generate $W'$ with these $n-k+r$ linearly independent vectors.\footnote{For details, see pg. 22 of \cite{quantum-book}} Since, the number of ways we can choose $r$ linearly independent vectors from $W$ is $(q^k-1)(q^k-q)\cdots(q^k-q^{r-1})$, the number of ways we can choose $n-k$ linearly independent vectors from $\mathbb{V}\setminus W$ is $(q^n-q^k)(q^n-q^{k+1})\cdots(q^n-q^{n-1})$ and the number of bases of an $n-k+r$ dimensional subspace is $(q^{n-k+r}-1)(q^{n-k+r}-q)\cdots(q^{n-k+r}-q^{n-k+r-1})$, the number of subspaces $W'$ with $dim(W')=n-k+r$ and $dim(W\cap W')=r$ is 
$$N_r=\dfrac{(q^k-1)(q^k-q)\cdots(q^k-q^{r-1})(q^n-q^k)(q^n-q^{k+1})\cdots(q^n-q^{n-1})}{(q^{n-k+r}-1)(q^{n-k+r}-q)\cdots(q^{n-k+r}-q^{n-k+r-1})}$$
Now, as $0\leq r \leq k-1$, 
$$deg(W)=\sum_{r=0}^{k-1}N_r.$$\qed
{\remark \label{max-degree-min-degree} It is clear from Theorem \ref{degree-theorem} that $deg(W)$ depends solely on its dimension and it is minimized if $dim(W)=1$ and maximized if $dim(W)=n-1$, i.e., $$\delta=q^{n-1} \mbox{ and }\Delta=\sum_{r=0}^{n-2}\left[ \begin{array}{c} n-1 \\ r \end{array}\right]_q \dfrac{q^{n-r-1}}{q^r+q^{r-1}+\cdots +q+1}.$$ }
{\corollary \label{eulerian-theorem} $\mathcal{G}(\mathbb{V})$ is Eulerian if and only if $q$ is even.}\\
\pf From Theorem \ref{degree-theorem}, 
$$N_r=\dfrac{\{(q^k-1)(q^k-q)\cdots(q^k-q^{r-1})\}\{(q^n-q^k)(q^n-q^{k+1})\cdots(q^n-q^{n-1})\}}{(q^{n-k+r}-1)(q^{n-k+r}-q)\cdots(q^{n-k+r}-q^{n-k+r-1})}$$
$$=\dfrac{\{q^{\{1+2+\cdots+(r-1)\}}(q^k-1)\cdots(q^{k-r+1}-1)\}\{q^{\{k+(k+1)+\cdots +(n-1)\}}(q^{n-k}-1)\cdots(q-1)\}}{q^{\{1+2+\cdots+ (n-k+r-1)\}}(q^{n-k+r}-1)(q^{n-k+r-1}-1)\cdots(q-1)}$$
$$=\dfrac{q^{r(r-1)/2}(q^k-1)\cdots(q^{k-r+1}-1)q^{\{n(n-1)-k(k-1)\}/2}(q^{n-k}-1)\cdots(q-1)}{q^{\{(n-k+r)(n-k+r-1)\}/2}(q^{n-k+r}-1)(q^{n-k+r-1}-1)\cdots(q-1)}$$
$$=\dfrac{q^{r(r-1)/2}q^{\{n(n-1)-k(k-1)\}/2}}{q^{\{(n-k+r)(n-k+r-1)\}/2}}\cdot (q^k-1)\cdots(q^{k-r+1}-1)(q^{n-k}-1)\cdots(q^{n-k+r+1}-1)$$
$$=q^{(n-k)(k-r)}\cdot (q^k-1)\cdots(q^{k-r+1}-1)(q^{n-k}-1)\cdots(q^{n-k+r+1}-1)$$

If $q$ is even, $N_r$ is even for all $r$ satisfying $0\leq r \leq k-1$. Hence, by Theorem \ref{degree-theorem}, for any subspace $W$ of $\mathbb{V}$ of dimension $k$, $deg(W)=\sum_{r=0}^{k-1}N_r$ is even. As all vertices of $\mathcal{G}(\mathbb{V})$ are of even degree, $\mathcal{G}(\mathbb{V})$ is Eulerian. On the other hand, if $q$ is odd, as the minimum degree $\delta=q^{n-1}$, $\mathcal{G}(\mathbb{V})$ is not Eulerian. \qed

{\corollary \label{edge-connectivity} Edge connectivity of $\Gamma(\mathbb{V})$ is $q^{n-1}$.}\\
\pf From \cite{plesnik}, as $\Gamma(\mathbb{V})$ is of diameter $2$ (by Theorem \ref{diameter}), its edge connectivity is equal to its minimum degree, i.e., $q^{n-1}$. \qed

{\theorem \label{clique-number-chromatic-number-theorem} Let $\mathbb{V}$ be an $n$-dimensional vector space over a finite field of order $q$. Then the clique number and chromatic number of $\mathcal{G}(\mathbb{V})$ are both equal to $1+q+\cdots+q^{n-1}$.}\\
\pf By Theorem \ref{maximal-clique}, $\mathcal{W}_{n-1}$ is a maximal clique and $$|\mathcal{W}_{n-1}|=\left[ \begin{array}{c} n \\ n-1 \end{array}\right]_q=\left[ \begin{array}{c} n \\ 1 \end{array}\right]_q=\dfrac{q^n-1}{q-1}=1+q+\cdots+q^{n-1}.$$ Now, by Lemma \ref{graph-homomorphism-lemma}, $\mathcal{G}(\mathbb{V})$ is $|\mathcal{W}_{n-1}|$-colourable, i.e., $\chi(\mathcal{G}(\mathbb{V}))\leq |\mathcal{W}_{n-1}|$. Moreover, as $\mathcal{W}_{n-1}$ is a maximal clique, $\omega(\mathcal{G}(\mathbb{V}))\geq |\mathcal{W}_{n-1}|$. As chromatic number of a graph is greater or equal to its clique number, we have the following inequality:
$$\omega(\mathcal{G}(\mathbb{V})) \leq \chi(\mathcal{G}(\mathbb{V}))\leq |\mathcal{W}_{n-1}| \leq \omega(\mathcal{G}(\mathbb{V})),$$ i.e., $$\omega(\mathcal{G}(\mathbb{V})) = \chi(\mathcal{G}(\mathbb{V}))= |\mathcal{W}_{n-1}| =1+q+\cdots+q^{n-1}.$$  \qed

{\remark \label{example-remark} Theorem \ref{clique-number-chromatic-number-theorem} shows that $\mathcal{G}(\mathbb{V})$ is weakly perfect. In Theorem \ref{perfect-theorem}, we establish a necessary and sufficient condition for $\mathcal{G}(\mathbb{V})$ to be prefect. In Figure \ref{diagram}, the $7$ inner vertices forms a maximum clique and the chromatic number is also $1+2+2^2=7$. Also, as $q=2$, by Corollary \ref{eulerian-theorem}, $\mathcal{G}(\mathbb{V})$ is Eulerian.}

{\theorem \label{perfect-theorem} $\mathcal{G}(\mathbb{V})$ is perfect if and only if $dim(\mathbb{V})=3$.}\\
\pf Let $n=dim(\mathbb{V})\geq 4$. Let $S=\{\alpha_1,\alpha_2,\alpha_3,\alpha_4\}$ be $4$ linearly independent vectors in $\mathbb{V}$ and let $S\cup T$ be the extension of $S$ to a basis of $\mathbb{V}$. Consider the following induced $5$-cycle in $\mathcal{G}(\mathbb{V})$: $W_1\sim W_2 \sim W_3 \sim W_4 \sim W_5 \sim W_1$ where $W_1=\langle \alpha_1,\alpha_2,T\rangle, W_2=\langle \alpha_3,\alpha_4,T\rangle, W_3=\langle \alpha_1+\alpha_3,\alpha_2,T\rangle, W_4=\langle \alpha_1,\alpha_4, T \rangle$ and $W_5=\langle \alpha_1+\alpha_3,\alpha_2+\alpha_4,T \rangle$. Thus, by Strong Perfect Graph Theorem, $\mathcal{G}(\mathbb{V})$ is not perfect.

Let $n=dim(\mathbb{V})=3$ and if possible, let there exists an induced odd cycle $W_1\sim W_2\sim W_3\sim \cdots \sim W_{2k+1} \sim W_1$ of  length greater than $3$ in $\mathcal{G}(\mathbb{V})$, where $W_i$'s are distinct proper non-trivial subspaces of $\mathbb{V}$.   Since, $W_i \sim W_{i+1}$, therefore $W_i + W_{i+1}=\mathbb{V}$, i.e., at least one of $W_i$ or $W_{i+1}$ is of dimension $2$. Without loss of generality, let $W_1$ be of dimension $2$. Now, there are two possibilities: $dim(W_2)=1$ or $dim(W_2)=2$.

If $dim(W_2)=1$, then as $W_2 \sim W_3$, $dim(W_3)$ must be $2$. However, in that case we have $W_1+W_3=\mathbb{V}$, i.e., $W_1 \sim W_3$, a contradiction. On the other hand, let $dim(W_2)=2$. Since, $W_1 \not\sim W_3$, i.e., $W_1+W_3\neq \mathbb{V}$, we have $dim(W_3)=1$. Again, as $W_2 \not\sim W_4$, i.e., $W_2+W_4\neq \mathbb{V}$, we have $dim(W_4)=1$. But, this implies $W_3+W_4 \neq \mathbb{V}$, i.e., $W_3 \not\sim W_4$, a contradiction. Thus, there does not exist any induced odd cycle of length greater than $3$ in $\mathcal{G}(\mathbb{V})$.

Let $n=dim(\mathbb{V})=3$ and if possible, let there exists an induced odd cycle $W_1\sim W_2\sim W_3\sim \cdots \sim W_{2k+1} \sim W_1$ of  length greater than $3$ in $\overline{\mathcal{G}(\mathbb{V})}$, where $W_i$'s are distinct proper non-trivial subspaces of $\mathbb{V}$. It is to be noted that in $\overline{\mathcal{G}(\mathbb{V})}$, $W_i \not\sim W_j$ implies $W_i+W_j=\mathbb{V}$. Thus, as in previous case, without loss of generality, let $W_1$ be of dimension $2$. Since, $W_2 \sim W_1$ in $\overline{\mathcal{G}(\mathbb{V})}$, $W_1+W_2\neq \mathbb{V}$, i.e.,  $W_2 \subset W_1$ and $dim(W_2)=1$.  Now, as $W_2 \not\sim W_4$ in $\overline{\mathcal{G}(\mathbb{V})}$, $W_2+W_4 = \mathbb{V}$, i.e., $dim(W_4)=2$. Also, as $W_3\sim W_4$, we have $W_3+W_4\neq \mathbb{V}$ and hence $dim(W_3)=1$. Similarly, as $W_5\sim W_4$, we have $W_5+W_4\neq \mathbb{V}$ and hence $dim(W_5)=1$. Again as $W_3\not\sim W_5$, we have $W_3+W_5= \mathbb{V}$. However this is impossible as $dim(W_3)=dim(W_5)=1$. Thus, there does not exist any induced odd cycle of length greater than $3$ in $\overline{\mathcal{G}(\mathbb{V})}$.

Thus, as neither $\mathcal{G}(\mathbb{V})$ nor its complement have any induced odd cycle of length greater than $3$, by Strong Perfect Graph Theorem, $\mathcal{G}(\mathbb{V})$ is perfect. Hence, the theorem. \qed

\section{Conclusion}
In this paper, we represent subspaces of a finite dimensional vector space as a graph and study various inter-relationships among $\mathcal{G}(\mathbb{V})$ as a graph and $\mathbb{V}$ as a vector space. The main goal of these discussions was to study the algebraic properties of $\mathbb{V}$ and graph theoretic properties of $\mathcal{G}(\mathbb{V})$, and to establish the equivalence between the corresponding graph and vector space isomorphisms. Apart from this, we also study basic properties e.g., completeness, connectedness, domination, chromaticity, maximal cliques, perfectness and maximal independent sets of $\mathcal{G}(\mathbb{V})$. As a topic of further research, one can look into the structure of automorphism groups, Hamiltonicity and independence number of $\mathcal{G}(\mathbb{V})$.

\section*{Acknowledgement}
The author is thankful to Sabyasachi Dutta, Anirban Bose and Sourav Tarafder for some fruitful discussions on the paper. The research is partially funded by NBHM Research Project Grant, (Sanction No. 2/48(10)/2013/ NBHM(R.P.)/R\&D II/695), Govt. of India.

\end{document}